\documentclass{amsart}
\usepackage{amssymb}
\usepackage{amsfonts}
\usepackage{amsbsy}

\usepackage{latexsym}
\usepackage[mathscr]{eucal}
\usepackage{verbatim}

\theoremstyle{plain}
\newtheorem{theorem}{Theorem}[section]
\newtheorem{corollary}[theorem]{Corollary}
\newtheorem{lemma}[theorem]{Lemma}
\newtheorem{proposition}[theorem]{Proposition}
\newtheorem{remark}[theorem]{Remark}
\newtheorem{definition}[theorem]{Definition}

\newcommand{\A}{\mathcal{A}}
\newcommand{\B}{\mathcal{B}}

\newcommand{\K}{\mathcal{K}}

\DeclareMathOperator{\TA} {\mathcal{T}(\A)}
\DeclareMathOperator{\Aff} {Aff(\TA)}

\DeclareMathOperator{\Ext} {\partial_{e}(\TA)}
\DeclareMathOperator{\her} {her}
\DeclareMathOperator{\tr} {Tr}
\DeclareMathOperator{\T} {T}
\DeclareMathOperator{\St} {\A\otimes \K}
\DeclareMathOperator{\M}{\mathcal M(\mathcal A \otimes \mathcal K)}

\def\sideremark#1{\ifvmode\leavevmode\fi\vadjust{\vbox to0pt{\vss
\hbox to 0pt{\hskip\hsize\hskip1em
\vbox{\hsize2cm\tiny\raggedright\pretolerance10000
\noindent#1\hfill}\hss}\vbox to8pt{\vfil}\vss}}}

\newcommand{\be}{\begin{equation}\label}
\newcommand{\ee}{\end{equation}}
\newcommand{\bq}{\begin{equation*}}
\newcommand{\eq}{\end{equation*}}
\newcommand{\ba}{\begin{align*}}
\newcommand{\ea}{\end{align*}}
\newcommand{\bp}{\begin{proof}}
\newcommand{\ep}{\end{proof}}
\newcommand{\bL}{\begin{lemma}\label}
\newcommand{\eL}{\end{lemma}}
\newcommand{\bP}{\begin{proposition}\label}
\newcommand{\eP}{\end{proposition}}
\newcommand{\bC}{\begin{corollary}\label}
\newcommand{\eC}{\end{corollary}}
\newcommand{\bT}{\begin{theorem}\label}
\newcommand{\eT}{\end{theorem}}
\newcommand{\bR}{\begin{remark}\label}
\newcommand{\eR}{\end{remark}}
\newcommand{\bD}{\begin{definition}\label}
\newcommand{\eD}{\end{definition}}

\numberwithin{equation}{section}
\author{Victor Kaftal}

\address{Department of Mathematics\\
University of Cincinnati\\
P. O. Box 210025\\
Cincinnati, OH\\
45221-0025\\
USA}

\email{victor.kaftal@uc.edu}

\author{P. W. Ng}

\address{Department of Mathematics\\
University of Louisiana\\
217 Maxim D. Doucet Hall\\
P.O. Box 41010\\
Lafayette, Louisiana\\
70504-1010\\
USA}

\email{png@louisiana.edu}

\author{Shuang Zhang}

\address{Department of Mathematics\\
University of Cincinnati\\
P.O. Box 210025\\
Cincinnati, OH\\
45221-0025\\
USA}

\email{shuang.zhang@uc.edu}

\date{4/30/2013}


\begin{document}
\title[Strict comparison of projections]{Strict comparison of projections and positive combinations of projections in certain multiplier algebras}
\begin{abstract} In this paper we investigate whether positive elements
in the multiplier algebras of certain finite C*-algebras can be written as finite linear combinations of projections with
positive coefficients (PCP). Our focus is on the category of
underlying C*-algebras that are  separable,
simple,  with real rank zero,  stable rank one, 
finitely many extreme traces, and strict comparison of projections by the traces. We  prove that the strict comparison of projections holds also in the
multiplier algebra $\M$. Based on this result and under the additional hypothesis   that $\M$ has real
rank zero, we characterize which  positive elements of $\M$ are of PCP.

\end{abstract}
\maketitle
\section{Introduction}\label{S:1}
In this article we focus on three closely related problems on
C*-algebras:
\begin{itemize}
\item [(A)] Which elements are (finite) sums of commutators?
\item [(B)] Is every element a (finite) linear combination of projections?
\item [(C)] Which positive elements are (finite) linear combination of projections with positive coefficients (called a positive combination of projections, or PCP for short)?
\end{itemize}
During the last several decades   much work have been done on these
problems  for various algebras, in particular for  (A) and (B). In
1954, Halmos proved that every bounded operator on an infinite
dimensional separable Hilbert space $H$ is a sum of two commutators.
Then in 1967 Fillmore \cite{Fillmore} found that every element of
$B(H)$ is a linear combination of 257 projections, eventually
reduced to 8 by Goldstein and Paszkiewicz \cite
{GoldsteinPaszkiewicz}. Fillmore also observed (\cite {Fillmore})
that  a positive compact operator of  infinite rank  cannot be a
positive combination of projections; Fong and Murphy proved (\cite
{FongMurphy}) that these are the only exceptions in $B(H)$.

Almost at the same time the same problems were investigated in von
Neumann algebras. In 1967 Pearcy and Topping
 \cite {PearcyToppingIdempotents} proved that every element in a properly infinite algebra is a sum of
 2 commutators and every self-adjoint element is a linear combination of 8 projections.
 In   1968   they proved \cite {PearcyTopping} that in every finite type I algebra,  self-adjoint elements in the kernel of
  the central trace are finite sums of commutators. The same result was proven for the more
  delicate  type II$_1$ von Neumann algebras  by Fack and De La Harpe
\cite {FackDelaHarpe}. Building on that,  Goldstein and Paszkiewicz
  showed \cite {GoldsteinPaszkiewicz}  that every
 element in a type II$_1$ von Neumann algebra is a linear combination of projections, but that every element in  a
 finite type I  algebra  is a linear combination of projections if and only if the center of the algebra is finite dimensional.

The investigation of (A) for C*-algebras started  in 1982 with the
work of Fack    who proved \cite {Fack} that self-adjoint elements
in properly infinite unital C*-algebras or in stable C*-algebras are
the sum of 5 self-commutators. Furthermore,   Fack  proved \cite
{Fack}  that in a simple unital AF algebra self-adjoint elements in
the kernel of all tracial states
 are sums of 7 commutators.  Thomsen   extended his work to a large class of AH algebras which
 includes the irrational rotation algebras and crossed products of Cantor minimal systems (\cite{ThomsenCommutator}).

An important extension of these results was obtained by  Marcoux
  who proved \cite {MarcouxSmallNumberCommutators}  that in simple
 C*-algebras of real rank  zero with a unique tracial state which gives the strict comparison of projections, every self-adjoint element
in the kernel of that trace is a sum of 2 self-commutators. In \cite
{MarcouxIndianaSpanProjections} he proved that  all commutators are
linear combinations of projections under mild conditions. This and
other work brought an affirmative answer to  (B) for all C*-algebras
in the following categories:
\begin{itemize}
\item Simple purely infinite C*-algebras;
\item AF-algebras with finitely many extremal tracial states;
\item AT-algebra with real rank zero and finitely many extremal tracial states.
\item Certain AH-algebras with real rank zero, bounded dimension growth, and finitely many extremal tracial states.
\end{itemize}

The first result  about (C) was the work of Fong and Murphy \cite
{FongMurphy} in $B(H)$. In recent years the  authors of the present
paper in \cite {KNZPISpan} made significant progress by proving that
every positive element in a purely infinite simple $\sigma$-unital
C*-algebra $\A$ or in its multiplier algebra $\mathcal M(\mathcal
A)$  is a positive combination of projections (a PCP).

For von Neumann algebras the answer to  (C) can be found in \cite
{KNZFiniteSumsVNA} mirrored the results in $B(H)$: a positive
element in a type II$_\infty$ factor is a PCP if either its range
projection is finite or the element does not belong to the ideal
generated by  finite projections (the Breuer ideal). The non factor
case is similar but   in terms of central ideals.

More complex is in C*-algebras of finite types. In \cite {KNZComm}
  under the assumption that $\A$ is simple, separable,
unital C*-algebras with real rank zero, stable rank one, strict
comparison of projections, and finitely many extremal tracial
states, we proved that  a positive element $a\in \A\otimes \K$  is a
PCP if and only if $\bar\tau(R_a)<\infty$ for every tracial state
$\tau$ on $\A$, where $\bar\tau$ denotes the extension of
$\tau\otimes \tr$ to a normal semifinite trace on $(\A\otimes
\K)^{**}$  and $R_a $  denotes the range projection of $a$ in
$(\A\otimes \K)^{**}$.

Aim of the present paper is to investigate Problem (C) for positive
elements of the multiplier algebra $\M$ for a   C*-algebra $\A$ of
finite type considered in \cite {KNZComm}. Of independent interest
we also consider (A) and (B) for in the corners of $\M$, as
necessary steps in investigating (C).

A key tool in  \cite {KNZComm} was to  the strict comparison of
projections in $\A$ and in $\A\otimes \K$.  In Section 3 we
establish a certain strict comparison of projections in $\M$ based
on the strict comparison of projections in $\A$. We will show that
strict comparison of projections holds also in $\M$ (Theorem \ref
{T:strict comparison in M}) in the sense that if $P$ and $Q$ are
projections in $\M$ with $Q\not \in \A\otimes \K$ and if
$\bar\tau(P)< \bar\tau(Q)$ for all tracial states $\tau$ on $\A$ for
which $\bar\tau(Q)< \infty$, then $P\precsim Q$.

It follows from the work of Fack and Marcoux     that every
self-adjoint element $T\in \M$
 is a sum of two self-commutators of elements of $\M$. In the case that $T\in P(\M)P$ for
some projection $P\in \M$,  it is natural to ask whether those
self-commutators can be chosen in $P(\M)P$.  However, a necessary
condition is that $\bar \tau(T)=0$ for every $\tau$ for which
$\bar\tau(P)<\infty$. One of our main results, Theorem \ref {T: sum
of comm} is that this condition is also sufficient. As a consequence
of this result we then obtain in Theorem \ref {T:lin comb}  that
every element $P(\M)P$ is a linear combination of projections in
$P(\M)P$ {\it with control on the coefficients} (see Section \ref
{S:5}). The control on the coefficients permits us to prove (\cite
{KNZPISpan} and \cite {KNZComm}) that every positive locally
invertible element in $P(\M)P$ is PCP (Corollary \ref {C:corner}).

In the case that $\M$ has real rank zero, by modeling the proof on
the proof of \cite [Theorem  6.1]{KNZComm}  we prove (Theorem \ref
{T:PCP}) that a positive element $T\in \M$ is a PCP if and only if
$\bar\tau(R_T)< \infty$ for every $\tau\in \TA$ for which $T\in
I_\tau$. Here $I_\tau$ is the closed ideal generated by the
projections of $\M$ with finite $\bar \tau$ values. This result is
the natural analog of the Fong and Murphy   result in
$B(H)$ in view of $B(H)= \mathcal M(\mathbb C\otimes \K)$.  Indeed  the usual trace $\tr$ on $B(H)$ is the extension $\bar \tau$ of the unique tracial state
$\tau$ on $\mathbb C$,  
$\K = I_\tau$,  and  $T\in B(H)$ is a PCP
if either $T\not \in \K$ or if $T\in \K$ with $\tr(R_T)<\infty$, i.e.,  $T$ has finite rank  (\cite {FongMurphy}). \\

The second author would like to acknowledge that his research is
partially supported by a travel grant from the Simons Foundation.

\section{Preliminaries}\label{S:2}
In this paper  $\A$  is always assumed to be  unital  separable
simple C*-algebra with real rank zero and  stable rank one unless
otherwise specified. We assume henceforth that $\A$ has some tracial
states and  denote by $\TA$ the tracial simplex of $\A$, that is the
collection of  tracial states on $\A$.  It is well known
that $\TA$ is a w*-compact convex set. Denote by $\Ext$ the extreme
boundary of $\TA$, that is the collection of extreme points of
$\TA$, or extreme traces for short.

Every tracial state $\tau\in \TA$ extends uniquely to the faithful, semifinite, normal trace $\tau\otimes \tr$ on $(\A\otimes \K)_+$. Notice that  $\tau\otimes \tr(p)> 0$ for any nonzero projection $p\in \A\otimes K$ and $\tau\otimes \tr(1\otimes e)=1$ for any rank-one projection $e\in \K$. Up to scalar multiples, all  semifinite, normal trace on $(\A\otimes \K)_+$ arise in this way. Thus without risk of confusion, we will identify $\TA$ with the collection of semifinite, normal traces  on $(\A\otimes \K)_+$, normalized by $\tau(1\otimes e)=1$ for any rank one projection $e\in \K$.

Throughout this paper we further assume that $\A$ has the property of strict comparison of projections property, namely that if $p$ and $q$ are projections in $\A$ and $\tau(p)< \tau(q)$ for all $\tau \in \TA$ then $p \prec q$ ($p$ is subequivalent to $q$ in the Murray-von Neumann sense).
As we have shown in \cite[Lemma 2.4]{KNZComm}, if strict comparison of projections holds for $\A$, it holds also for $\A\otimes \K$.

Next for every $\tau\in \TA$ denote by $\tilde \tau$ the natural extension of $\tau$ to $\A^{**}$, then $\tilde \tau$ is a normal tracial state on $\A^{**}$ and hence $\bar\tau := \tilde\tau \otimes \tr$ is a normal semifinite (not necessarily faithful) trace (a tracial weight) on $\A^{**}\otimes B(H)$ which we can identify with $(\A\otimes \K)^{**}$. Thus $\bar \tau$ is also  a semifinite trace on the multiplier algebra $\M_+ $.
Notice that as remarked in  \cite[5.3]{KNZComm}, by the work of F. Combes   \cite [Proposition 4.1 and Proposition 4.4] {Combes} and Ortega, Rordam, and Thiel
\cite[Proposition 5.2]{OrtegaRordamThiel}
$\tau\otimes \tr $ has a unique extension to a lower semicontinuous semifinite  trace $\bar \tau$  on the enveloping von Neumann algebra $(\A\otimes \K)^{**}_+$ and hence this extension is $\bar\tau$.

Recall that every  open projection  $P\in (\A\otimes \K)^{**}$, and in particular  every projection $P\in  \M$, has a decomposition  $P=\bigoplus_1^\infty p_j$ into a series of strictly  converging  projections $p_j\in \A\otimes \K$. Thus
$\bar\tau(P) = \sum_1^\infty \tau(p_j)$ for all $\tau\in \TA$.

One of the goals of this paper will be to show (see Theorem \ref {T:strict comparison in M}) that under the additional hypothesis that $\Ext$ is finite, $\M$ has a form of strict comparison of projections with respect to the traces $\{\bar \tau\mid \tau\in \TA\}$ as described below.

\bD{D:strict} If $P$ and $Q$ are projections in $\M$ we say that $P$ is tracially dominated by $Q$ if
$$ \bar\tau(P)< \bar\tau(Q)\quad \text{for all } \tau\in \TA \text { for which } \bar\tau(Q)< \infty.
$$
\eD
Under the convention that ``$\infty < \infty$", $P$ is tracially dominated by $Q$ if $ \bar\tau(P)< \bar\tau(Q)$ for all $ \tau\in \TA$.

Notice that if  $\TA$ has only a finite number of extremal traces, say $\Ext = \{\tau_j\}_1^n$, then $ \bar\tau(P)< \bar\tau(Q)$ for all $ \tau\in \TA$  if and only if $\bar\tau_j(P)< \bar\tau_j(Q)$ for all $1\le j\le n$.

Strict comparison of projections within $\A\otimes \K$ means  that if $p$ and $q$ are projections in $\A\otimes \K$ and  $p$ is tracially dominated by $q$, then $p\prec q$.

The same cannot hold without further conditions in $\M$ and the obstruction arises from the ideal structure of $\M$. Indeed, if $P\precsim Q$ then $P$ must belong to the same ideal as $Q$. But it will be easy to show that for every nonzero projection $q\in \A\otimes \K$, which necessarily generates $\A\otimes \K$ as an ideal, there is a projection $P\not \in  \A\otimes \K$  that is tracially dominated by $q$. Thus before further considering strict comparison of projections in $\M$ we need to recall  some facts about  ideals in $\M$.

First of all, there exists a minimal ideal $I_{min}$ that properly contains $\A\otimes K$ (\cite [Theorem 1.7]{LinZhang}).  The projections in $I_{min}$ are the characterized  by the following property.

\bD{D:thin} \cite[Definition 2.1]{LinZhang} A sequence of projections $p_j\in \A\otimes \K$ is called an $\ell^1$ -sequence if for every projection $0\ne r\in \A\otimes \K$ there is an  $N\in \mathbb N$ such that
$$[p_n]+[p_{n+1}]\cdots [p_m]\le [r] \quad \forall ~m>n\ge N.$$ A projection  $P\in \M$ is called {\it thin} if it has a decomposition $P=\bigoplus_1^\infty p_j$  into a  strictly converging sum of an $\ell^1$ -sequence.
\eD

We collect here the following known results.
\bP{P:thin}\cite {LinZhang}
\item [(i)] If $P\in \M$ is thin, then for
any  decomposition  $P=\bigoplus_1^\infty p_j$ into a series of strictly converging projections   $p_j\in \A\otimes \K$, the sequence $\{p_j\}$ is $\ell^1$. \item [(ii)] Finite sums of thin projections are thin and projections majorized by thin projections are thin.
\item[(iii)] Every thin projection $P\in \M\setminus \A\otimes K$ generates the ideal $I_{min}$ and every projection in $I_{min}$ is thin.
\eP

Furthermore, recall that  for every projection $p\in \A\otimes \K$, the evaluation map $\hat p: \TA\ni \tau\to \tau(p)$ is affine and continuous. Let $\Aff$ denote the Banach space of real-valued affine continuous functions on $\TA$. Recall that under the evaluation map, $K_o(\A)$  is dense in $\Aff$ (\cite [Theorem 6.9.3] {BlackadarBook}) and  the dimension semigroup $D(\A\otimes K)\setminus \{0\}$ is dense in $\Aff_{++}$, the collection of strictly positive affine continuous functions on $\TA$ (see \cite [Remark 2.7]{KNZComm}.)

For every projection $P\in \M$, the evaluation map $\hat P: \TA\ni \tau\to \bar\tau(P)$ is affine and lower semicontinuous.  Indeed if $P=  \bigoplus_1^\infty p_j$ is any strictly convergent decomposition of $P$ with $p_i\in \A\otimes K$, then $\bar \tau(P)= \sum_1^\infty \tau(p_i)$, that is, $\hat P=\sum_1^\infty\widehat {p_i}$, as the pointwise sum of a series of positive continuous affine functions, is affine lower semicontinuous.

\bL{L:cont}
 Let $P\in \M$. Then $P\in I_{min}$ if and only if the evaluation map $\hat P: \TA\ni \tau\to \bar\tau(P)$ is continuous. Furthermore, if $P=  \bigoplus_1^\infty p_j\in I_{min}$ is any strictly convergent decomposition of $P$ with $p_i\in \A\otimes K$, then the series $ \sum _1^\infty \hat p_j$ converges uniformly.
\eL
\bp
Assume first that $P\in  I_{min}$. Hence by Proposition \ref {P:thin} (iii), $P=  \bigoplus_1^\infty p_j$ is the sum of an $\ell^1$ sequence and it will be enough to prove that the series $ \sum _1^\infty \hat p_j$ converges uniformly.   For every $\epsilon >0$, by \cite [Theorem 6.9.3] {BlackadarBook} choose a nonzero projection $r\in \A\otimes \K$ with $\tau(r)< \epsilon$ for all $\tau\in \TA$. Then there is an $N\in \mathbb N$ such that for all $m>n\ge N$ ,$ \bigoplus_n^m p_j \precsim r$ and hence $\sum _n^m \hat p_j(\tau) < \tau(r) < \epsilon$. Thus the sequence of partial sums of the continuous functions $\hat p_j$ is uniformly Cauchy and hence its limit  $\hat P$ is continuous.

Conversely, assume that $\hat P$ is continuous. Then by Dini's theorem the sequence of partial sums  $\sum _1^n  \hat p_j$ increases uniformly to $\hat P$. Let $r\in \A\otimes \K$ be a nonzero projection. Since $(\hat r)(\tau)>0$ for all $\tau$, $\hat r$ is continuous, and $\TA$ is compact, then there is some $\alpha>0$ such that $(\hat r)(\tau)\ge \alpha >0$ for all $\tau$. But then there is an $N\in \mathbb N$  such that for all $m>n\ge N$,
$$\tau( \bigoplus_n^m p_j) = \bigoplus_n^m \hat p_j(\tau) < \alpha\le  \tau(r)\quad \forall ~\tau\in \TA.$$
By the strict comparison of projections in $\A\otimes\K$ it follows that
$\bigoplus_n^m p_j\precsim r$ for all $m>n\ge N$ and hence $P$ is thin. By Proposition \ref {P:thin} (iii), $P\in I_{min}$.
 \ep

Whenever $\TA$ does not reduce to a singleton (a unique trace) there are more proper ideals between $\A\otimes \K$ and  $\M$.

\bD{D:ideals}
For every $\tau\in \TA$, denote by $I_\tau$ the closed ideal of $\M$ generated by all the projections  $P\in \M$ with $\bar\tau(P)<\infty$.
\eD
Recall from \cite [Theorem 2.3] {ZhangRiesz} that every projection $P\in  I_\tau$ satisfies the condition $\bar\tau(P)<\infty$.

Furthermore, as a consequence of results in \cite {ZhangRiesz}
and \cite{ZhangMatricial} $$I_\tau=\overline{ \{X\in \M\mid  \bar\tau(X^*X)< \infty\}}=\text{span}\overline{ \{X\in \M_+\mid  \bar\tau(X)< \infty\}} $$
where the closures are in norm.

In the case when $\Ext$ is finite, say $\Ext= \{\tau_j\}_1^n$, then there are precisely $2^n-1$ ideals distinct from $\A\otimes \K$, which are obtained by all the possible intersections of the maximal ideals $I_{\tau_j}$ for $1\le j\le n$. This was obtained in \cite {LinIdeals}  for AF algebras and in \cite{RordamIdeals} for the general case.

In the case when $\Ext$ is infinite the ideal structure of $\M$ is considerably more complex. Some results in particular when $\TA$ is a Bauer simplex, are obtained in \cite {PereraIdeals}

In section 5 we will make use of the following simple result

\bL{L:stuff} Let $T\in \M_+$ and let $\tau\in \TA$.
Then  $T\in I_\tau$ if and only if $\bar\tau(\chi_{(\delta, \|T\|]}(T))<\infty$ for every $\delta>0$.
\eL
\bp
Assume first that $\bar\tau(\chi_{(\delta, \|T\|]}(T))<\infty$ for every $\delta>0$. Let $$f_\delta(t)=\begin{cases}t-\delta&\text{if } t> \delta\\ 0 &\text{if } t \le  \delta\end{cases} $$ and let $T_\delta:=f_\delta(T)$. Then $R_{T_\delta}= \chi_{(\delta, \|T\|]}(T)$ and hence  $\bar\tau(T_\delta)\le \|T_\delta\| \bar\tau(R_{T_\delta})< \infty$. Thus $T_\delta\in I_\tau$ and since  $\|T-T_\delta\|\le \delta$ for every $\delta>0 $ and $I_\tau$ is closed, it then follows that   $T\in I_\tau$.

Now assume that $T\in I_\tau$. By the definition of $I_\tau$, for every $\delta>0$ there is a $B\in (I_\tau)_+$ with $\tau(B)<\infty$ and  $\|T-B\|<\frac{\delta}{4}$. By basic results about Cuntz subequivalence (see for example \cite [Lemma 2.2.] {RordamIdeals}) it follows that $(T-\frac{\delta}{2} 1)_+\precsim B$. But then there is an $x\in \M$ such $(T-\frac{\delta}{2} 1)_+=x^*Bx$. Hence
$$\bar\tau((T-\frac{\delta}{2} 1)_+)\le \bar\tau(x^*Bx)\le \bar\tau(B^{1/2}x^*xB^{1/2})\le \|x\|^2\bar\tau(B)< \infty.$$
But then
$$ (T-\frac{\delta}{2} 1)_+= (T-\frac{\delta}{2}1)\chi_{(\frac{\delta}{2}, \|T\|]}(T)\ge (T-\frac{\delta}{2} 1)\chi_{(\delta, \|T\|]}(T)\ge \frac{\delta}{2} \chi_{(\delta, \|T\|]}(T)$$ and hence $\bar\tau (\chi_{(\delta, \|T\|]}(T))< \infty$.
\ep

\section{Strict comparison of projections in $\M$}\label{S:3}

We start by considering strict comparison of of projections belonging to the ideal $I_{min}$.
\bP{P:str comp I_o}
Let $P,Q\in I_{min}$ be projections for which $P$ is tracially dominated by $Q$ and $Q\not \in \A\otimes \K$. Then $P\prec Q$.
\eP

\bp
Let $P=\bigoplus_1^\infty p_j$,  $Q=\bigoplus_1^\infty q_j$ with $p_j, ~q_j$ projections in $\A\otimes \K$ and the series converging strictly, and furthermore $q_j\ne 0$ for all $j$ by the assumption that  $Q\not \in \A\otimes \K$.  Then by Proposition \ref {P:thin} and Lemma \ref {L:cont} these sequences are $\ell^1$ and  the series of continuous affine functions $\sum_1^\infty \widehat{p_j}$ and $\sum_1^\infty \widehat{q_j}$ converge uniformly.

By a routine compactness argument, we can find an index $N$ such that
$$\hat P(\tau) < \sum_1^{N} \widehat{q_j}(\tau) \quad \forall ~ \tau\in \TA.$$
If only finitely many projections $p_j$ are nonzero, then $P\in \A\otimes \K$. But then by the strict comparison of projections in $\A\otimes \K$ it follows that
$P\prec \bigoplus_1^{N}q_j\le Q$ and we are done. Thus assume that all projections $p_j\ne 0$.

Since the series of continuous functions $\sum_{1}^\infty \widehat{p_j}$ converges uniformly and  $q_{N+1}\ne 0$, we can find an index $m_1$ such that
$$\|\sum_{j=m_1+1}^\infty \widehat{p_j}\|_\infty < \min \widehat{q_{N+1}}.$$
Thus for all $\tau\in \TA$ we have
\be{e:step1}
\sum_1^{m_1} \widehat{p_j}(\tau)~ < ~\hat P(\tau)~ < ~\sum_1^{N} \widehat{q_j}(\tau)\quad\text{and}\quad
\sum_{m_1+1}^\infty \widehat{p_j}(\tau)~ < ~\widehat{q_{N+1
}}(\tau).
\ee
Since $\lim _m\|\sum_{j=m}^\infty \widehat{p_j}\|_\infty=0$ and  $\min \widehat{q_{N+k}}>0$ for all $k$, choose an increasing sequence of indices $m_k$ such that for all $k\ge 1$ and all $\tau\in \TA$
$$
\sum_{m_k+1}^{m_{k+1}} \widehat{p_j}(\tau)~ < ~\widehat{q_{N+k}}(\tau).
$$
By the strict comparison of projections in $\A\otimes \K$, we obtain  that
 $\bigoplus_1^{m_1} p_j\prec \bigoplus_1^{N} q_j$ conjugated by a partial isometry $V_0\in \M$ and
 $\bigoplus_{m_k+1}^{m_{k+1}} p_j\prec q_{N+k}$, conjugated by a partial isometry $V_k\in \M$ for all $k\ge 1$. By the strict convergence of the series $P=\bigoplus_1^\infty p_j$ and $Q=\bigoplus_1^\infty q_j$, it follows that also the series $\sum_{k=0}^\infty V_k$ converges strictly. Thus its sum $W:=\sum_{k=0}^\infty V_k$ is a partial isometry in $\M$, $WPW^*\le Q$, and hence $P\prec Q$.
 \ep
\

To consider comparison of projections not in $I_{min}$ we will need to further assume that $\Ext$ is finite. In that case, by the complete characterization of ideals in $\M$ (see \cite {RordamIdeals}) it follows that      $$I_{min}= \underset{\tau\in \TA}{\bigcap} I_{\tau}= \underset{\tau\in \Ext}{\bigcap} I_{\tau}.$$

\bT{T:strict comparison in M}
Assume that $\TA$ has  finite extreme boundary $\Ext=\{\tau_j\}_1^n$. Then strict comparison of projections holds for $\M$ in the sense that if  $Q\not \in \A\otimes \K$ and $\bar\tau(P)< \bar\tau(Q)$ for all $\tau\in \TA$ for which $\bar\tau(Q)< \infty$ then $ P\precsim Q.$
\eT
\bp
Let $\Ext=\{\tau_j\}_1^n$ be the extremal boundary of $\TA$ and let $$S: =\{1\le j\le n \mid  \bar\tau_j(Q)
< \infty\}.$$ If $S=\emptyset$, i.e., $\bar\tau_j(Q)=\infty$ for all $j$, it follows that $\bar\tau(Q)=\infty$ for all $\tau\in \TA$. Then $Q\sim 1$  by \cite [2.4 or 3.6] {ZhangRiesz}  and thus $P\precsim Q$. If $S=  \{1,2, \cdots, n\}$, i.e., $\bar\tau(Q)< \infty$ for all $\tau\in \Ext$, then by the remark preceding this theorem,  $Q\in I_{min}$ and hence $P\prec Q$ by Proposition \ref {P:str comp I_o}. Thus assume henceforth that $\emptyset\subsetneq S\subsetneq  \{1,2, \cdots, n\}$.

Let $\alpha:= \min_{j\in S} \big(\bar\tau_j(Q)- \bar\tau_j(P)\big)$. Then $\alpha >0$. Let $P=\bigoplus_1^\infty p_j$,  $Q=\bigoplus_1^\infty q_j$ with $p_j, ~q_j$ projections in $\A\otimes \K$,  the series converging strictly, and $q_j\ne 0$ for all $j$.
Since the series $\sum_i^\infty \tau_j(q_i)$ converges for all $j\in S$, we can find an integer $n_0$ such that
$\tau_j(q_{n_0})< \alpha$ for all $j\in S$. Then
$$\sum_1^\infty \tau_j(p_i)=\bar\tau_j(P)\le \bar\tau_j(Q)-\alpha< \bar\tau_j(Q)-\tau_j(q_{n_0})=  \sum_1^{\infty}\tau_j(q_i)- \tau_j(q_{n_o})\quad \forall ~j\in S.
$$
Thus there is an integer $n'\ge  n_o$ for which
\be{e: 1} \sum_1^\infty \tau_j(p_i)<  \sum_1^{n'}\tau_j(q_i)- \tau_j(q_{n_o})\quad \forall ~j\in S.
\ee
By hypothesis, $q_i\ne 0$ for all $i$ and in particular as all the traces are faithful,  $\tau_j(q_{n_0})>0$ for all $j\in \{1, 2, \cdots, n\}$.
Since the series $\sum_i^\infty \tau_j(p_i) $ converge for all $j\in S$, we can then find $m_1$ such that
\be{e:2}\sum_{m_1+1}^\infty \tau_j(p_i) <\tau_j(q_{n_0})\quad\forall ~j\in S.\ee
By (\ref {e: 1}),
$$  \sum_1^{m_1}  \tau_j(p_i)<  \sum_1^{n'}\tau_j(q_i)- \tau_j(q_{n_o})\quad \forall ~j\in S.$$ By using the divergence of the series
$\sum_1^{\infty} \tau_j(q_i)$ for $j\not \in S$, we can also find an $n_1> n'$ such that
\be{e:3} \sum_1^{m_1}  \tau_j(p_i)<\sum_1^{n_1-1}\tau_j(q_i)- \tau_j(q_{n_o}) \quad \forall~j\in \{1, 2, \cdots, n\}.
\ee
This concludes the initial step. Now choose $m_2>m_1$ and  $n_2> n_1$ such that
\begin{alignat}{2}
\sum_{m_2+1}^\infty \tau_j(p_i)&< \tau_j(q_{n_1})\quad \forall~j\in S && (\text{the left series converges})\notag\\
\sum_{m_1+1}^{m_2}\tau_j(p_i)&< \sum_{i=n_1+1}^{n_2-1}\tau_j(q_i) +\tau_j(q_{n_o})  ~ \forall~j\in \{1, 2, \cdots, n\} ~ &&(\text{the right series diverges}).\label{e:4}
\end{alignat}
Iterating the construction, we can find an increasing sequence of indices $m_k$ and $n_k$ such that
\be{e:5}
\sum_{m_k+1}^{m_{k+1}} \tau_j(p_i)< \sum_{ n_k+1}^{n_{k+1}-1}  \tau_j(q_i) +\tau_j(q_{n_{k-1}}) \quad \forall~j\in \{1, 2, \cdots, n\}.
\ee
Then by the strict comparison of projections in $\A\otimes \K$ we have for all $k$
\begin{alignat*}{2}
\bigoplus_1^{m_1}p_i&\prec \bigoplus_1^{n_1-1}q_i - q_{n_o} &&(\text{by (\ref {e:3})})\\
\bigoplus_{m_k+1}^{m_{k+1}}p_i&\prec \bigoplus_{ n_k+1}^{n_{k+1}-1}q_i +q_{n_{k-1}} \quad &&(\text{by (\ref {e:5})})
\end{alignat*}
Reasoning as in the proof of Proposition \ref {P:str comp I_o}  we can construct a partial isometry in $\M$ to conjugate $P$ to a subprojections of $Q$, thus obtaining that $P\precsim Q$.

\ep

\bP{P:Dim map onto}
Assume that $\TA$ has  finite extreme boundary $\Ext=\{\tau_j\}_1^n$. Then for every n-tuple of $\alpha_j\in (0, \infty]$, there exist a projection $P\in \M\setminus \A\otimes \K$ such that $\bar\tau_j(P)=\alpha_j$ for all $1\le j\le n$.
\eP
\bp
Let $S:= \{ 1\le j\le n \mid  \alpha_j<\infty\}$. If $S=\emptyset$, it is enough to choose $P=1$. To simplify notations, assume  that $\emptyset \ne S\ne \{1,2, \cdots, n\}$,  the proof in the case when $S=\{1,2, \cdots, n\}$ being identical.

Let $1=\bigoplus_1^\infty E_i$ be a strictly converging decomposition of the identity into projections $E_i\sim1$. Recall that for any infinite collection of nonzero projections $P_j\le E_j$, the sum $\bigoplus_1^\infty P_i$ converges in the strict topology to a projection $P\in \M\setminus \St$.

Recall also that $\Aff$ is isomorphic to $\mathbb R^n$ and that the dimension semigroup $\mathrm{D}(\A\otimes \K)$ is dense in $\Aff_{++}$, the collection of strictly positive continuous affine functions on $\TA$.
Thus there is a projection $p_1\in \A\otimes\K$ such that
$$
\begin {cases}\alpha_j -1 < \tau_j(p_1)<  \alpha _j& j\in S\\
1< \tau_j(p_1)< 2 &j\not \in S.
\end{cases}
$$
Since $p_1\prec  1\sim E_1$,  we can choose $p_1\le E_1$.
Next, we find a projection $p_2 \in \A\otimes\K$ with $p_2\le E_2$ and
$$
\begin {cases} \alpha_j -\tau_j(p_1)-\frac{1}{2}< \tau_j(p_2)<  \alpha _j- \tau_j(p_1)& j\in S\\
1< \tau_j(p_2)< 2 & j\not \in S.
\end{cases}
$$
Iterating, we find a sequence of  projection $p_i\in \A\otimes \K$, with $p_i\le E_i$ for which
$$\begin{cases} \alpha_j- \frac{1}{m}< \sum_{i=1}^m\tau_j(p_i)< \alpha_j & j\in S\\
 \sum_{i=1}^m\tau_j(p_i)> m& j\not \in S.
  \end{cases}
$$
But then $\bigoplus_{i=1}^\infty p_i$ converges to a projection $P\in \M\setminus \A\otimes \K$ and
$\bar\tau_j(P)=\alpha_j$ for every $1\le j\le n$.

\ep

By combining Proposition \ref {P:Dim map onto} with Theorem \ref {T:strict comparison in M} we thus obtain:

\bC{C: comp} Assume that $\TA$ has  finite extreme boundary $\Ext=\{\tau_j\}_1^n$.   For every projection $Q\in \M\setminus \A\otimes \K$ and every  n-tuple of $\alpha_j\in (0, \infty]$ with $\alpha_j< \bar\tau_j(Q)$ for all $j$ for which $\bar\tau_j(Q)< \infty$,  there is a projection $P$ in $\M\setminus \A\otimes \K$ such that $P\prec Q$ and $\bar\tau_j(P)=\alpha_j$ for all $1\le j\le n$.

\eC
\section{Sums of commutators in the ideals of $\M$} \label{S:4}

Fack proved in \cite [Theorem 2.1] {Fack} that if a unital algebra $\B$ contains two mutually orthogonal projections equivalent to the identity, then every selfadjoint element $b$ is the sum of five selfcommutators,
\be{e:Fackbounds} b= \sum_{i=1}^5 [x_i, x_i^*] \quad \text{with}\quad \|x_j\|\le \frac{3}{2}\|b\|^{1/2} \quad\text{for } 1\le j\le 5.\ee
The bound $ \|x_1\|\le \frac{3}{2}\|b\|^{1/2}$ for the element obtained in the first step of the proof is implicit in \cite [Lemma 1.2] {Fack}, while the bounds $\|x_i\|\le\|b\|^{1/2}$ for the remaining elements can be seen from the proofs of  \cite [Lemma 1.3, 1.4, 2.3] {Fack}.

Since the identity of $\M$ can be decomposed into the sum of two  mutually orthogonal projections equivalent to the identity, it follows that every element of $\M$ is the sum of 10 selfcommutators.

In particular, every element in a corner  $P\M P$ for some projection $P\in \M$ is a sum of commutators of elements of $\M$. These elements don't necessarily belong to $P\M P$. Indeed  if $\bar \tau(P)< \infty$ for some $\tau\in \TA$, then $\bar \tau$  is a finite trace on $P\M P$ and hence vanishes on all the commutators of elements of  $P\M P$.  Thus for $T\in P\M P$ to be a sum of commutators of elements of  $P\M P$ it is necessary that $\bar \tau(T)=0$ for all $\tau\in \TA$ for which $\bar\tau (P)< \infty$. We shall prove that the condition is also sufficient under the additional hypothesis that  $\Ext$ is finite.
Based on the work in \cite {LinRR0AHEmbedding}, \cite {MarcouxSmallNumberCommutators}, and \cite {ThomsenCommutator},  we obtained in a previous paper:

 \bL{lem:MarcouxPropositionReplacement} \cite[Lemma 3.3]{KNZComm} Let $\B$ be a unital separable simple C*-algebra
of real rank zero, stable rank one, and  strict comparison of
projections.  Let $b \in \B$  be a selfadjoint element, let
$\eta>0$, and assume that  $ |\tau(b)|\le \eta$ for all $ \tau \in
T(\B)$. Then for every $\epsilon > 0$ there exist elements $v_1, v_2, v_3,
v_4 \in \B$ such that $\| v_i \| \le \sqrt{2} \| b \|^{1/2}$ for $1
\le i \le 4$ and
$$\| b - \sum_{i=1}^{4} [v_i, v_i^*] \| < \eta +\epsilon.$$
\eL

We start with the following extension of this lemma to the corners  of $\A\otimes \K$ by projections of $\M$.
\bL {L: corner of stab A}
Assume that $\TA$ has  finite extreme boundary $\Ext=\{\tau_j\}_1^n$,  let $P$ be a nonzero projection in $\M$,  and set $$S: =\{1\le j\le n  \mid \bar\tau_j(P)< \infty\}\quad \text{and} \quad  \alpha := \underset{j\in S}{\min }~ \bar\tau_j(P). $$ Let $a \in P(\A\otimes \K) P$  be a selfadjoint element, let
$\eta>0$, and assume that  $ |\tau_j(a)|\le \eta$ for all $ j \in
S$. Then for every $\epsilon>0$  there exist  elements $v_1, v_2, v_3,
v_4, v_5 \in P(\A\otimes\K) P$ such that
$\|v_i\|\le \sqrt 2 \|a\|^{1/2}\quad \text{for } 1\le i \le 5$ and
$$\|a- \sum_{j=1}^5[v_i,v_i^*]\| < \frac{\eta}{\alpha}+\epsilon.$$
\eL

\bp
The case when $S=\emptyset$, namely when $\bar\tau_j(P)=\infty $  for all $1\le j\le n$ and hence $\bar\tau(P)=\infty $ for all $\tau\in \TA$, is immediate because then $P\sim 1$ and without loss of generality,  $P=1$. But then by Fack's \cite [Theorem  1.1]{Fack}, $a$ is the sum of five selfcommutators. The bounds on the norms are implicit in Fack's proof.

We leave to the reader the case when $S=\{1, 2, \cdots  n\}$ which is similar but somewhat simpler than the general case. Thus assume  that $\emptyset \ne S\ne \{1, 2, \cdots  n\}$. To simplify notation, assume furthermore that $\|a\|=1$.

Choose $0< \beta < \alpha$ so that
\be{e: beta} \frac{\eta}{\beta}< \frac{\eta}{\alpha}+ \frac{\epsilon}{5}.
\ee

Decompose $a$ into its positive and negative parts $a=a_+-a_-$. Since $\A\otimes K$ has real rank zero, by \cite[Lemma 2.3]{KNZ Mult}  we can approximate  from underneath $a_+$ (resp., $a_-)$ with a positive finite spectrum element, that is,  find   $\sum_{i=1}^{n_1} \lambda_i p_i\le a_+$ with mutually orthogonal projections $p_i\in \A\otimes K$ and $\lambda_i>0$, (resp., $\sum_{k=1}^{n_2} \mu _kq_k\le a_-$ with mutually orthogonal projections $q_k\in \A\otimes K$ and $\mu_k>0$) so that
\be {e:def b} b:= a_+- \sum_i^{n_1} \lambda_i p_i -a_-+ \sum_{k=1}^{n_2} \mu _kq_k\ee
has norm
$ \|b\|\le \frac {\epsilon\beta }{ 5\max _{j\in S}\bar\tau_j(P)}$.
In particular,
\be{e:norm b}
\|b\|< \frac {\epsilon }{ 5}.
\ee
Set
$$
a':= \sum_i^{n_1} \lambda_i p_i-\sum_{k=1}^{n_2} \mu _kq_k.
$$
Since $R_b\le P$,  for all $j\in S$ we have
$ |\tau_j(b)|\le \frac{\epsilon\beta}{5}.$
Since $a'= a- b$ and hence $|\tau_j(a')|\le |\tau_j(a)| + |\tau_j(b)|$, we also have
\be{e:tau(a')}
|\tau_j(a')|\le \eta+  \frac{\epsilon\beta}{5}\quad \forall~j\in S.
\ee
Choose an integer $m \ge \frac{5}{\epsilon}$. By the density of the dimension semigroup $\mathrm{D}(\A\otimes \K)$  in the collection  $\Aff_{++}$ of strictly positive continuous affine functions on $\TA$, choose projections in $\{p_i'', q_k''\}$ in $\A\otimes K$ with
\begin{align*}&
\begin{cases}\tau_j(p_i) - \frac{\beta}{(n_1+n_2)(m+1)}< \tau_j(p_i'')< \tau_j(p_i) & j\in S\\
\tau_j(p_i'')< \min \{\frac{\beta}{(n_1+n_2)(m+1)},  \tau_j(p_i)\} &j\not\in S
\end{cases}\\
&\begin{cases}\tau_j(q_k) - \frac{\beta}{(n_1+n_2)(m+1)} < \tau_j(q_k'')< \tau_j(q_k) & j\in S\\
\tau_j(q_k'')< \min \{\frac{\beta}{(n_1+n_2)(m+1)}, \tau_j(q_k)\} &j\not\in S.
\end{cases}
\end{align*}
Then by using strict comparison of projections in $\A\otimes K$ (see Corollary \ref {C: comp}), find projections $p_i'\sim p_i''$, $q_k'\sim q_k''$ in $\A\otimes K$ with
$p_i'\le p_i$, $q_k'\le q_k$ for all $i, k$. Set
\begin{align*}
r&:= \sum_{i=1}^{n_1} (p_j- p'_i)+ \sum_{k=1}^{n_2}(q_k-q'_k)\\
r'&:= \sum_{i=1}^{n_1}  p'_i+ \sum_{k=1}^{n_2}q'_k\\
c&:=  \sum_{i=1}^{n_1} \lambda_i (p_i-p'_i)-\sum_{k=1}^{n_2} \mu _k(q_k-q'_k)\\
c'&:= \sum_{i=1}^{n_1} \lambda_i p'_i-\sum_{k=1}^{n_2} \mu _kq'_k.
\end{align*}

Notice that  the projections $\{p_i, q_k\}$, and hence the projections $\{p_i', q_k'\}$. are all mutually orthogonal and are majorized by $P$, hence $r$ and $r'$ are projections in $P(\A\otimes K)P$,  and $c$ and $c'$ are selfadjoint elements of $P(\A\otimes\K)P$ with range projections $R_c=r$ and $R_{c'}=r'$ respectively. Then
\begin{align*}
 a'&= c+c' \\
\tau_j(r)&<  \frac{\beta}{m+1}   \text{ for } ~ j \in S\\
\tau_j(r')&< \frac{\beta}{m+1} \text{ for }  j\not\in S.  \\
 \intertext{ Since $R_c=r$ and $\|c\|\le \|a\|=1$, $R_{c'}=r'$ and $\|c'\|\le \|a\|=1$, we have }
|\tau_j(c)|&\le \|c\|\tau_j(r)\le \frac{\beta}{m+1}< \frac{\epsilon \beta}{5}\quad\text{for } j \in S \\
|\tau_j(c')|&\le \|c'\|\tau_j(r')\le \frac{\beta}{m+1}< \frac{\epsilon \beta}{5}  \quad\text{for } j\not \in S.
\end{align*}

By (\ref {e:tau(a')}) and the above inequalities we obtain
$$
|\tau_j(c')|\le \begin{cases}  |\tau_j(a')|+ | \tau_j(c)|< \eta + \frac{2\epsilon\beta}{5} & j\in S\\
 \frac{\epsilon \beta}{5} &j\not \in S
\end{cases} \qquad \le \eta + \frac{2\epsilon\beta}{5}\quad \forall~j
$$
and hence
\be{e:tau(c')2}  |\tau(c']|< \eta + \frac{2\epsilon\beta}{5}\quad \forall~\tau\in \TA.\ee

Since  by construction  $r'\in P(\A\otimes \K)P$ and hence $r'\ne P$, it follows that  $\tau_j(r')< \bar\tau_j(P)$ for all $j$. Invoking again the density of $\mathrm{D}(\A\otimes \K)$  in   $\Aff_{++}$, choose a projection  $s\in\A\otimes \K$ with
$$
\max(\tau_j(r'), \beta)< \tau_j(s) < \begin{cases}\bar \tau_j(P)&j\in S\\
\max(\tau_j(r'), \beta)+1&j\not\in S\end{cases}
$$
Using strict comparison of projections in $\A\otimes \K$ and in $\M$, it is now routine to show that $s$ can be chosen so that $r'\le s\le P$. By construction, $\tau(s) > \beta$ for all $\tau\in \TA$.

Now  $c'$ belongs to $s(\A\otimes \K) s$ which is a unital separable simple C*-algebra of real rank zero, stable rank one, and with strict comparison of projections.
Every tracial state $\tilde \tau\in \T(s(\A\otimes \K) s)$ is the restriction and rescaling of a trace in $\tau\in \TA$, i.e.,
$ \tilde \tau(c')= \frac{\tau(c')}{\tau(s)}$. But then  for every $ \tilde \tau\in \T(s(\A\otimes \K) s)$ we have by (\ref {e: beta}) and (\ref {e:tau(c')2})
$$
| \tilde \tau(c')|< \frac{|\tau(c')|}{\beta} < \frac{ \eta +\frac{2\epsilon\beta }{5}}{\beta}= \frac{\eta}{\beta} +  \frac{2\epsilon}{5}<  \frac{\eta}{\alpha} +  \frac{3\epsilon}{5}
$$
Thus by   \cite[Lemma 3.3]{KNZComm} (see Lemma \ref {lem:MarcouxPropositionReplacement}), we can find elements  $v_1, v_2, v_3, v_4$ in $s (\A\otimes \K) s\subset  P(\A\otimes \K) P$ with
 \begin{align}& \|v_i\|\le \sqrt 2\|c'\|^{1/2}\le\sqrt 2\|a\|^{1/2}=\sqrt 2 \quad\text{for all } ~ 1\le i\le 4\label{e:norm bounds}\\
&\| c'- \sum_{i=1}^4 [v_i,v_i^*] \|<   \frac{\eta}{\alpha}+ \frac{3\epsilon}{5}. \label{e:norm2}
  \end{align}

Now we consider $c$ which has range projection $r\in \A\otimes\K$. Since
$$\tau_j(r) < \begin{cases} \frac{\beta}{m+1} < \frac{\bar\tau_j(P)}{m+1} & j\in S\\  \infty =  \frac{\bar\tau_j(P)}{m+1} & j\not\in S,
\end{cases}
$$
by Theorem \ref {T:strict comparison in M}
 we can find $m$ mutually orthogonal  subprojections $\{r_j\}_1^m$ of $P-r$ with $r_j\sim r$. Since $\A\otimes \K $ is an ideal of $\M$, all the projections $r_j$ are also in $\A\otimes \K$. Let $r_j=v_jrv_j^*$ for some partial isometries $v_j\in \A\otimes \K $. Then we can identify $e:=c- \frac{1}{m}\sum_{i=1}^m v_icv_i^*$ with the diagonal element of $\mathbb M_{m+1}(r(\A\otimes \K) r)$ given by the matrix

$$\begin{pmatrix} c&0&\cdots& 0\\
0 & -\frac{1}{m}c&\cdots&0\\
&&\ddots&&\\
0&\cdots&0 &
-\frac{1}{m}c
\end{pmatrix}$$

From \cite [Fack, Lemma 1.3]{Fack}, there is an element $v_5\in \mathbb M_{m+1}(r(\A\otimes \K) r)$, which in turns we can identify with  an element of $(\sum_{i=0}^m r_i)(\A\otimes \K )(\sum_{i=0}^mr_i) \subset P(\A\otimes \K) P$, where $r_o:=r$, such that
$e= [v_5,v_5^*]$ and $\|v_5\|\le \|c\|^{1/2}\le  \|a\|^{1/2}=1$.
But
\be{e:norm3}
\|c- [v_5,v_5^*]\|=  \frac{\|c\|}{m}\le \frac{\|a\|}{m}\le \frac{\epsilon}{5}.
\ee
We thus have
$$a- \sum_{i=1}^5 [v_i,v_i^*]= b + c'- \sum_{i=1}^4 [v_i,v_i^*] + c- [v_5,v_5^*],$$
and hence from (\ref {e:norm b}), (\ref{e:norm2} ), and (\ref {e:norm3})
$$\|a- \sum_{i=1}^5 [v_i,v_i^*]\|\le \frac{\eta}{\alpha} +\epsilon.$$
\ep

Now we extend Lemma \ref {L: corner of stab A} to corners of $\M$.

\bL{L:corner of M}
Assume that $\TA$ has  finite extreme boundary $\Ext=\{\tau_j\}_1^n$,  let $P$ be a nonzero projection in $\M$,  and set $$S: =\{  1\le j\le n \mid  \bar\tau_j(P)< \infty\}. $$ Let $T \in P\M P$  be a selfadjoint element and assume that  $ \tau_j(T)=0$ for all $ j \in
S$. Then for every $\epsilon>0$  there are 10  elements $V_1, V_2, \cdots V_{10}\in P\M P$ such that
$\|V_i\|\le 4\sqrt 2 \|T\|^{1/2}~\, \text{for } 1\le i \le 10$ and
$$\|T- \sum_{j=1}^{10}[V_i,V_i^*]\| < \epsilon.$$
\eL
\bp
Let $\pi: P\M P \to P\M P/P(\A\otimes\K)P$ be the canonical quotient map. By \cite[Theorem 3.6]{KuchPer} and for a more general case, \cite [Theorem A and 4.5]{KuchNgPer}, the corona algebra $\M/\A\otimes \K$ is purely infinite and hence so is $P\M P/P(\A\otimes\K)P$. For completeness purpose, we show how this fact follows easily from Theorem \ref{T:strict comparison in M}.

By \cite {ZhangMatricial} one can write $P=Q\oplus Q'$, where $Q\sim Q'$ are projections of $\M$.
Then $\bar \tau (Q)=\bar \tau (Q') = \frac 12 \bar \tau (P)$ for
all $\tau \in \TA$. Write $ P =\bigoplus_1^\infty p_j$ as a strictly
convergent sum of projections of $\A\otimes \K$. Then $\bar \tau (P) =\sum_{j=1}^\infty \tau (p_j)$
for all $\tau \in \TA$.

Let $s_n = \bigoplus_1^n p_j$ and choose $n_o$ so that for all $n\ge n_o$, $$\bar \tau_j (P-s_n) \begin{cases} < \bar \tau_j
(Q)<\infty& \forall j\in S\\  = \bar \tau_j (Q)
=\infty& \forall j\not \in S\
\end{cases}.$$ It follows from Theorem \ref{T:strict comparison
in M} that $P-s_n\precsim Q$. Since $s_n\in \A
\otimes \K$, it follows that $\pi (P)= \pi (P-s_n)\prec \pi(Q)$ as wanted.

Then by Fack's \cite [Theorem 2.1]{{Fack}} and (\ref {e:Fackbounds}), there are five  elements $$\tilde V_i\in P\M P/P(\A\otimes\K)P$$ such that
$$
\pi(T) =\sum_{i=1}^5 [\tilde V_i,\tilde V_i^*]\quad \text{and} \quad
\|\tilde V_i\|\le \frac{3}{2} \|\pi(T)\|^{1/2}\le \frac{3}{2} \|T\|^{1/2}\quad \text{for } ~1\le i\le5.
$$
Now choose liftings $V_i\in P\M P$ of $\tilde V_i$ such that $\|V_i\|\le \frac{2}{\sqrt 3}\|\tilde V_i\|$, hence $\|V_i\|\le \sqrt 3 \|T\|^{1/2}$ for $1\le i\le 5$.
Let
\be{e:def a}a:= T- \sum_{i=1}^5 [V_i, V_i^*]. \ee
Then $a=a^*\in P(\A\otimes K)P$. Notice that $$ -V_i^*V_i\le V_iV_i^*-V_i^*V_i \le  V_iV_i^*$$ and hence $\| [V_i, V_i^*]\|\le \|V_i\|^2.$ Thus
$$
\|a\|\le \|T\|+ \sum_{i=1}^5 \|V_i\|^2\le 16\|T\|.
$$
Furthermore, $\tau_j(V_iV_i^*)\le \|V_i\|^2\bar\tau_j(P)< \infty$ for every $j\in S$. Then
$$\tau_j(a)=\tau_j(T) - \sum_{i=1}^5 \tau_j([V_i, V_i^*]) =0\quad\forall~j\in S.$$
By Lemma \ref {L: corner of stab A}, the selfadjoint element $a$ can be approximated by the sum of five selfcommutators  of elements $V_j\in P(\A\otimes \K)P$ with
$$\|V_j\|\le \sqrt 2 \|a\|^{1/2}\le 4\sqrt 2\|T\|^{1/2} \quad\text{for } ~6\le i\le10$$ and
$$\|a-\sum_{j=6}^{10} [V_i, V_i^*]\|< \epsilon.$$
This combined with (\ref {e:def a}) concludes the proof.

\ep

Notice that the bounds are of course far from sharp. By same proof we could replace $4\sqrt2$ with any number strictly larger than $3.5\sqrt 2$.
 \bL {L:halving} Assume that $\TA$ has  finite extreme boundary $\Ext=\{\tau_j\}_1^n$,  let $P$ be a nonzero projection in $\M\setminus \A\otimes \K$,  and set $$S: =\{1\le j\le n \mid \bar\tau_j(P)< \infty\}. $$ Then there are three sequences of projections $P_n, Q_n, R_n$ in $P\M P$ with the properties
\item [(i)] $P_1+Q_1+R_1=P.$
\item [(ii)]The projections $\{R_n\}_1^\infty$ are mutually orthogonal.
\item [(iii)] $P_n+Q_n=R_{n-1}$ for all $n\ge 2$.
\item[(iv)] $\begin{cases}P_1\sim Q_1\prec R_1& n=1\\
P_n \sim Q_n\sim R_n &n\ge 2.
\end{cases}$
\item [(v)] $\tau_j(R_n)=\tau_j(P_n)= \tau_j(Q_n)=\infty$ for every $j\not \in S $ and every $n$.
\eL
\bp
By using the fact that every projection  in $\M\setminus  \A\otimes \K$ can be halved, i.e., decomposed into the sum of two orthogonal equivalent projections \cite [Theorem 1.1]{ZhangMatricial}  it is routine to find a sequence $\{R_n\}$ of mutually orthogonal subprojections of $P$ for which
$$
2[R_n]= \begin{cases}[P]& n=1\\
[R_{n-1}]& n\ge 2
\end{cases}
$$
Then by halving we find projections $P_n\sim Q_n$ so that
$$
P_n+Q_n= \begin{cases}P-R_1& n=1\\
R_{n-1} & n\ge 2
\end{cases}
$$
For $n=1$ we have $2[P_1]= [P-R_1]=[R_1]$ hence $P_1\sim Q_1\prec R_1$, while
for $n\ge 2$ we have $2[P_n]= [R_{n-1}]= 2[R_n]$, hence $P_n \sim Q_n\sim R_n $. (v) is now immediate since $\bar\tau_j(P)=\infty$ for all $j\not\in S$.
\ep

\bT {T: sum of comm}
Assume that $\TA$ has  finite extreme boundary $\Ext=\{\tau_j\}_1^n$,  let $P$ be a nonzero projection in $\M\setminus \A\otimes \K$,  and set $$S: =\{1\le j\le n \mid \bar\tau_j(P)< \infty\}. $$
Let $T=T^*\in P\M P$, and assume that $\tau_j(T)=0$ for all $j\in S$. Then there are two elements $V_i\in P(\M) P$ such that
 $T= \sum_{j=1}^{2}[V_i,V_i^*].$
 Furthermore, the elements $V_j$ can be chosen such that $\|V_j\|\le c\|T\|^{1/2}$ where $c $ is a constant that does not depend  on $T$, $P$, or the C*-algebra $\A$.
\eT
\bp

By adapting Fack's  proof of
\cite[Theorem 3.1]{Fack} and its modification by Thomsen \cite [Theorem 1.8]{ThomsenCommutator}, Marcoux has shown in
 \cite [Lemma 3.9] {MarcouxSmallNumberCommutators}  that if $\B$ is a simple, unital C*-algebra with real rank zero, strict comparison of projections and a unique tracial state, then every  selfadjont element in the kernel of the trace is the sum 8 selfcommutators. The two key elements of his proof are  \cite[Lemma 3.7]{Fack}, which state the existence of a sequence of projections satisfying  conditions (i)-(iii) and a weaker form of (iv) of Lemma \ref  {L:halving},  and  \cite[Proposition 3.6]{Fack}, which states that every selfadjont element in the kernel of the trace can be approximated arbitrarily well by the sum of 2 selfcommutators (with control on the norms of the operators). If we replace the latter result by Lemma \ref  {L:corner of M} which provides an approximation by the sum of 10 selfcommutators (with control on the norms of the operators), we see that  Marcoux's  proof  holds  in our setting and shows that  every selfadjont element $T\in P\M P$ in the kernel of the extremal traces $\{\tau_j\}_{j\in S}$  is the sum of 16 selfcommutators in $P\M P$ (with control on the norms of the operators).

Moreover, the number of selfcommutators can be reduced to two as in  \cite [Theorem 3.10] {MarcouxSmallNumberCommutators}. The norms of the elements forming the selfcommutators are bounded by a constant multiple $c$ of $\|T\|^{1/2}$,  but the estimates of $c$ are very far from sharp as discussed in \cite [Remark 5.3]{MarcouxSmallNumberCommutators} (see also  \cite [Proof of Theorem 3.4, Remark 3.5]{KNZComm}.
\ep

\section{Linear combination of projections in $\M$} \label{S:5}

Marcoux has shown in \cite [Theorem 3.8]{MarcouxIndianaSpanProjections} that in every unital C*-algebra that contains three projections $P_1+P_1+P_3=1$ with $P_i\precsim 1-P_i$ for $1\le i\le 3$,  every commutator $[x,y]$  is a linear combination of  84 projections. Furthermore, in \cite [5.1]{MarcouxSmallNumberCommutators}, Marcoux notices that the coefficients in this linear combination can be bounded by $8\|x\|\,\|y\|$.

Since the identity of $1\in \M$ can be decomposed into the sum of three mutually orthogonal projections $P_i\sim1$, and since, as remarked at the beginning of Section \ref {S:3}, every element $T\in \M$ is the sum of 10 commutators (with control on the norms), it follows immediately that every element of $T\in \M$ is a linear combination of 840 projections with control on the norms of the coefficients. However the results in the previous sections permit us to obtain a stronger result, namely that if $T$ is in a corner $P\M P$ of $\M$ then $T$ is a linear combination  of projections belonging to the same corner.

\bT{T:lin comb} There are constants $N$  and $M $ such that if
$\A$ is a unital separable simple C*-algebra
of real rank zero, stable rank one, strict comparison of
projections and has finitely many extreme tracial states and $P\in \M$ is a projection,  then every element $T\in P\M P$ is a linear combination of projections $T= \sum_{j=1}^N\lambda_jP_j$ with $\lambda_j\in \mathbb C$ and $P_j\le P$ projections, with $|\lambda_j|  < M$ for all $j$.
\eT
\bp
In the case when $P\in \A \otimes \K$ and hence $P\M P= P(\A \otimes \K)P$, the result follows from \cite [Theorem 4.4] {KNZComm}. Thus assume henceforth that $P\not \in \A \otimes \K$.

Let $\Ext = \{\tau_j\}_1^n$ and let $S: =\{j\in \mathbb N\mid  1\le j\le n, ~ \bar\tau_j(P)< \infty\}.$  Assume also that $T\ge 0$ and $0\ne\| T\|< 1$.

By Proposition \ref {P:Dim map onto} there is a projection $Q\in \M\setminus \A\otimes \K$ such that
$$\bar \tau(Q) = \begin{cases} \bar\tau(T) &j\in S\\
1&j\not \in S
\end{cases}
$$
Since $ \tau_j(T)< \bar\tau_j(P)$ for all $j\in S$, it follows by  Theorem \ref {T:strict comparison in M} that $Q\prec P$, so assume without loss of generality that $Q\le P$. Let $B:=T-Q$. Then \linebreak$B=B^*\in P\M$ and $\bar\tau_j(B)=0$ for all $j\in S$.  By Theorem \ref {T: sum of comm}, $B$ is the sum of two selfcommutators (with control on the norms of the elements), and each is the linear combination of  84 projections by \cite [Theorem 3.8]{MarcouxIndianaSpanProjections} with control of the coefficients (see remarks preceding this theorem), thus $T$ is a linear combination of 169 projections in $ P\M P$, also with control on the coefficients. As a consequence, every $T \in P\M P$ is a linear combination of 676 projections $P\M(\A\otimes\K) P$, also with control on the coefficients.

\ep
We say that an algebra $\B$ is the linear span of its projections {\it with control on the coefficients}  if it has a constant $V$ such that for every $b\in \B$ there are $n$ scalars  $\lambda_j\in \mathbb C$ and projections $p_j\in \B$ such that
\begin{itemize}
\item[(i)]  $b=\sum_1^n\lambda_jp_j$;
\item [(ii)] $\sum_1^n | \lambda_j |< V \|b\|$.
\end{itemize}

Thus Theorem \ref {T:lin comb} states that if the extremal boundary $\Ext$ is finite, then every hereditary subalgebra $P\M P$ of $\M$ is the linear span of its projections with control on the coefficients.

\section{Positive combination of projections in $\M$} \label{S:6}
Now we start investigating linear combinations of projections with positive coefficients, ({\it positive combinations of projections}, or PCP for short). We are interested in the question of which, necessarily positive, elements are PCP.

We obtained in \cite [Proposition 2.7]{KNZPISpan}  extending a $B(H)$ result by Fong and Murphy \cite {FongMurphy}, that if $\B$ is a unital C*-algebra  that is the span of its projections  with control on the coefficients and if PCPs are norm dense in $\B_+$, then every positive invertible element of $\B$ is a PCP.

Even in the case when $\M$ does not have real rank zero,  PCPs are norm dense  in $\M^+$ by \cite [Theorem 1.1] {Zhang}.  The same holds for  all the corners  $P\M P$ for projections $P\in \M$. Thus combining Theorem \ref {T:lin comb} and \cite [Proposition 2.7]{KNZPISpan} we obtain the following result.

\bC{C:corner} Assume that $\TA$ has  finite extreme boundary. Then for every projection $P\in \M$, every positive invertible element of $P\M P$ is a PCP.
\eC

As in \cite {KNZPISpan} and \cite {KNZComm}, the key tool for constructing PCP decompositions in our setting is given by the following result, which is an immediate consequence of  Corollary \ref {C:corner} and \cite [Lemma 2.9] {KNZPISpan}:

\bL{L: 2X2} Assume that $\TA$ has  finite extreme boundary.
 Let  $P, Q$ be projections in $\M$ with $PQ=0$, $Q\precsim P$ and let  $B=QB=BQ$ be
 a positive element of $\M$. Then for every  scalar $\alpha
>\|B\|$, the positive element $T:=\alpha P\oplus B$ is a PCP.
\eL

The next step is to prove  that if $\M$ has real rank zero, every positive element in a corner $P\M P$ that has sufficiently large range with respect to $P$ is also PCP (see statement  below). The proof is modeled on one of \cite[Lemma 6.4]{KNZComm} but with some substantial differences, so for clarity and completeness sake, we present a proof here.

\bL{L:cornerPCP}   Assume that $\TA$ has  finite extreme boundary and that $\M$ has real rank zero. Let  $P\in \M$ be a projection and $T\in P\M_+ P$ satisfy the conditions
\item [(i)] for every $\tau\in \TA$,  $T\in I_\tau$ if and only if $\bar\tau(P)< \infty;$
\item [(ii)]
$\begin{cases}\bar\tau(R_T)>\frac{1}{2}\bar\tau(P) & \text{if } \bar\tau(P)< \infty\\ \bar\tau(R_T)= \infty & \text{if }  \bar\tau(P)= \infty.\end{cases}$\\
Then  $T$ is a PCP.
\eL
 \begin{proof}
The case when $P\in \A\otimes \K$ is covered by \cite[Lemma 6.4]{KNZComm}, thus assume henceforth that $P\not \in \A\otimes \K$.  To simplify notations, assume  that $\|T\|=1$.
Let $\Ext= \{\tau_i\}_1^n$ and let $S: =\{1\le j\le n  \mid \bar\tau_j(P)< \infty\}. $ Assume that $\emptyset \ne S \ne \{1\le j\le n\}$, leaving to the reader the simpler cases when $S=\emptyset$ and when $S= \{1\le j\le n\}$.

Consider $T$  as an element of $P\M P$ and denote by $\chi(T)$ its spectral measure with values being projections in  $(P\M P)^{**}$.

By the w*-lower semicontinuity of each $\bar\tau_i$ and the w*-continuity of the restriction of $\bar\tau_i$ to $P\M P$ for each $ i\in S$, we have
\begin{align*}\lim_{\lambda\to0+}\bar\tau_i(\chi_{( \lambda, 1]}(T))&
=\bar\tau_i(R_T)    &&\forall  1\le i\le n \hspace{2cm}  \\
\lim_{\lambda\to0+}\bar\tau_i(\chi_{(0, \lambda)}(T))&=0  &&\forall  i\in S.\end{align*}
For every $i\not \in S$, by hypothesis $T\not\in I_{\tau_i}$, hence by Lemma \ref{L:stuff}   there is a $\delta_i>0$ such that $\bar \tau_i(\chi_{(\delta_i, 1]}(T))=\infty$. Let $\gamma_6:=\underset{i\not\in S}{\min}\, \delta_i$. Then
\be{e:ineq1}
  \bar\tau_i(\chi_{(\gamma_6,1]}(T)) \ge \chi_{(\delta_i,1]}(T)= \infty \quad  \forall~i\not\in S.
\end{equation}

For every $i\in S,$
$$\bar\tau_i(\chi_{\{0\}}(T))=\bar\tau_i(P)- \bar\tau_i(R_T)< \bar\tau_i(R_T).
$$ Thus we can find $0< \gamma_4 < \gamma_6$ such
that
\be{e:ineq2}
 \bar\tau_i(\chi_{[0, \gamma_4)}(T)) < \bar\tau_i(\chi_{(\gamma_6,1]}(T))  \quad \forall~i\in S.
\end{equation}

Now choose numbers $\gamma_1, \gamma_2, \gamma_3, \gamma_5 $ so that $0 < \gamma_1< \gamma_2< \gamma_3< \gamma_4< \gamma_5 < \gamma_6< 1$.
Let $f : [0, 1] \rightarrow [0, 1]$ be the  continuous
function defined by
$$
f(t) =
\begin{cases}
t & t\in [0,1]\setminus [\gamma_1, \gamma_3]\\
\gamma_1 &  t \in[\gamma_1, \gamma_ 2]\\
\text{linear }  & t \in[\gamma_2, \gamma _3].
\end{cases}
$$
Since $RR(P\M P)=0$, by Brown's interpolation property \cite{BrownIMP}, there exist
projections $S, R, Q \in  P\M P $ such that
\begin{alignat*}{5}
&\chi_{[0, \gamma_1]}(T)&~\le~ &S&~\le~& \chi_{[0, \gamma_2)}(T)\\
&\chi_{[0, \gamma_3]}(T)&~\le~ &Q&~\le~ &\chi_{[0, \gamma_4)}(T)\\
&\chi_{[\gamma_6, 1]}(T)&~\le~ &R&~\le~& \chi_{(\gamma_5, 1]}(T).
\end{alignat*}
Then
$$\bar\tau_i(Q)\le \bar\tau_i(\chi_{[0, \gamma_4)}(T))\le \bar\tau_i(\chi_{[ \gamma_6, 1]}(T))\le \bar\tau_i(R)\quad \forall ~i$$
and by (\ref {e:ineq2}), the inequality $\bar\tau_i(Q)\le \bar\tau_i(R)$ is strict for $i\in S$ while   $\bar\tau_i(R)= \infty $ for $i\not\in S$ by  (\ref {e:ineq1}). Then by Theorem \ref {T:strict comparison in M}  we obtain that $ Q\precsim R.$

Since $S-\chi_{[0, \gamma_1]}(T)\le \chi_{(\gamma_1, \gamma_2)}(T)$ and the function $f(t)$ is constant on the interval $[\gamma_1, \gamma_2]$, it follows that $S-\chi_{[0, \gamma_1]}(T)$ and hence $S$  commute with $f(T)$. Define:
\begin{align*}
T_1 &: = f(T)-f(T)S-\gamma_4 R\\
B~&:=T - f(T) + f(T) S\\
T_2&:=B+ \gamma_4 R.
\end{align*}

By a simple computation,
$$T_1\ge  \min\{ \gamma_1, \gamma_5-\gamma_4\}(P-S)$$ and
$R_{T_1}= P-S$.   Thus $T_1$ is positive and invertible in $(P-S)\M(P-S)$ and hence a PCP by Corollary \ref{C:corner}.

Since $B=QBQ\ge 0$, $QR=0$,  $ Q\precsim R$ and $\|B\|  \le  \gamma_2< \gamma_4,$   $T_2$ is a PCP by Lemma \ref{L: 2X2}. Since $T= T_1+T_2$  this concludes the proof.

\end{proof}

 \bT{T:PCP}
 Assume that $\TA$ has  finite extreme boundary and that $\M$ has real rank zero. Then $T\in \M_+$ is PCP if and only if either $T$ is full (i.e., belongs to no proper ideal of $\M$), or $\bar\tau(R_T)< \infty$ for every $\tau\in \TA$ for which $T\in I_\tau$.
 \eT
 \bp
 Assume that $T$ is PCP, namely $T =\sum_{j=1}^n\lambda_jP_j$ with scalars $\lambda_j>0$ and projections $P_j\in \M$. If  $T\in I_\tau$ for some $\tau\in \TA$, then for every $1\le j\le n$ it follows that $P_j\le \frac{1}{\lambda_j}T$ and hence $P_j\in  I_\tau$. But then $\bar\tau(P_j)< \infty$. Since $R_T= \bigvee_1^n P_j$, it follows by standard properties of traces on von Neumann algebras  that $$\bar\tau(R_T)\le \sum_{j=1}^n\bar\tau(P_j)< \infty.$$
 Assume now that $\bar\tau(R_T)< \infty$ for every $\tau\in \TA$ for which $T\in I_\tau$. Let $\Ext= \{\tau_i\}_1^n$ and let $S: =\{ 1\le j\le n \mid T\in I_{\tau_j}\}.$

 In the case when $S=\emptyset$,  $\bar\tau(R_T)=\infty$  for all $\tau\in \TA$, hence the result is given by Lemma \ref {L:cornerPCP} applied to $P=1$. Thus  assume that $S\ne \emptyset.$

By Proposition \ref {P:Dim map onto}, there is a projection $P\in \M\setminus \A\otimes \K$ such that $$\begin{cases}\bar\tau_j(R_T)< \bar\tau_j(P)< 2\bar\tau_j(R_T)< \infty & j\in S\\ \bar\tau_j(P)= \infty & j\not \in S.\end{cases}$$

Reasoning as in the proof of \cite [Lemma 6.3]{KNZComm} and using the strict comparison of projections in $\M$ (Theorem \ref  {T:strict comparison in M}), we can find a partial isometry $W\in \M^{**}$ such that $WW^*=R_T$, $W^*W\le P$ and such that the map $\Phi(X):=W^*XW$ is a *-isomorphism between $\her T = \her R_T$ and $\her (W^*W)\subset P\M P$. Since $\bar\tau(R_{\Phi(T)})=\bar\tau(R_T)$ for all $\tau\in \TA$, by Lemma \ref {L:cornerPCP}, $\Phi(T)$ is PCP in $P\M P$ and hence in $\her (W^*W)$. But then $T$ is PCP in $\her T $, which completes the proof.

 \ep


\begin{thebibliography}{99}
\bibitem{BlackadarBook} B. Blackadar, K-theory for operator algebras, second
edition,
Cambridge University Press, New York, 1998.


\bibitem{BrownIMP} L. G. Brown, \textit{Interpolation by projections
in C*-algebras of real rank zero,} J. Operator Theory, \textbf{26}
(1991) no. 2, 383--387.

\bibitem{Combes} F. Combes, \textit{Poids sur une C*-alg\'{e}bre,} J. Math. Pures Appl., IX. \textbf{47}  (1968) 57--100.

\bibitem{CowardElliottIvanescu} K. T. Coward, G. A. Elliott and
C. Ivanescu, \textit{The Cuntz semigroup as an invariant for
C*-algebras,} J. Reine Angew. Math., \textbf{623} (2008)
161--193.

\bibitem{CuntzPedersen}  J. Cuntz and G. K. Pedersen,
\textit{Equivalences and traces on C*-algebras,} J. Funct. Anal.,
\textbf{33} (1979) 135--164.



\bibitem{Fack} T. Fack, \textit{Finite sums of commutators in C*-algebras,}
Ann. Inst. Fourier, Grenoble, \textbf{32} (1982) 129--137.

\bibitem{FackDelaHarpe}  T. Fack and P. De la Harpe, \textit{Sommes de
commutateurs dans les algebres de von Neumann finies continues,} Ann. Inst.
Fourier, Grenoble, \textbf{30} (1980) 49--73.

\bibitem{Fillmore}  P. A. Fillmore, \textit{Sums of operators with square zero,} Acta Sci. Math. (Szeged), \textbf{28} (1967) 285--288.



\bibitem{FongMurphy} C. K. Fong and G. J. Murphy,
\textit{Averages of projections,} J. Operator Theory, \textbf{13}
(1985) no. 2, 219--225.

\bibitem{GoldsteinPaszkiewicz}  S. Goldstein and A. Paszkiewicz,
\textit{Linear combinations of projections in von Neumann algebras,}
Proc. Amer. Math. Soc., \textbf{116} (1992) no. 1, 175--183.


\bibitem{Halmos} P.R. Halmos, \textit{Commutators of operators. II} Amer. J. Math \textbf{76}, (1954), 191-198.


\bibitem{KNZPISpan} V. Kaftal, P. W. Ng and S. Zhang,
\textit{Positive combinations and sums of projections in purely infinite
simple C*-algebras and their multiplier algebras,} Proceedings AMS \textbf{139}, (2011), no 8, 2735-2746.

\bibitem{KNZChina}  \textit{Positive combinations of projections in von Neumann algebras and purely infinite simple C*-algebras}. Science China Mathematics, \textbf{54},(2011), no. 2,  1-10

\bibitem{KNZ Mult} V. Kaftal, P. W. Ng and S. Zhang,
\textit{Projection decomposition in multiplier algebras,} (2012). Math Ann
 \textbf {352}, (2012), no 3,  543-566.

\bibitem{KNZFiniteSumsVNA} H. Halpern, V. Kaftal, P. W. Ng and
S. Zhang, \textit{Finite sums of projections in von Neumann algebras,}
Transactions AMS \textbf{365}, (2013), 2409-2445

\bibitem{KNZTorsion}  V. Kaftal, P. W. Ng and S. Zhang,     \textit{Finite sums of projections in purely infinite simple C*-algebras with torsion K$_0$}, Proceedings AMS, \textbf{140}, (2012),  No 9, 3219-3227.

\bibitem{KNZComm} V. Kaftal, P. W. Ng and S. Zhang, \textit{Commutators and linear spans of projections in certain simple real C*-algebras}. Preprint (2012)  ArXiv:1208.1949

\bibitem{KuchPer} D. Kucerovsky and F. Perera,\textit{
Purely infinite corona algebras of simple C*-algebras with real
rank zero}
J. Operator Theory, \textbf{65}, (2011), no. 1, 131-144.

\bibitem{KuchNgPer} D. Kucerovsky, P. W. Ng and F. Perera,
\textit{Purely infinite corona algebras of simple C*-algebras},
Math. Ann. \textbf{346}, (2010) no. 1, 23--40.

\bibitem{LinIdeals}H. Lin, \textit{Ideals of multiplier algebras of simple AF C*-algebras.} \textbf{104}, (1988) Proc. Amer.Math. Soc.,
239--244.

\bibitem{LinRR0AHEmbedding} H. Lin, \textit{Embedding an AH-algebra into a
simple C*-algebra with prescribed KK-data,} K-theory, \textbf{24}
(2001) no. 2, 135--156.

\bibitem{LinCuntzSemigroup}  H. Lin, \textit{Cuntz semigroups of C*-algebras
of stable rank one and projective Hilbert modules,} (2010) Preprint,
A copy is available at
http://arxiv.org/pdf/1001.4558.

\bibitem{LinZhang} H. Lin, and S. Zhang, \textit{Certain simple C*-algebras with nonzero real rank whose corona algebras have real rank zero} (1992), Houston J. Math, No 1, 57-71.

\bibitem{MarcouxIndianaSpanProjections} L. W. Marcoux,
\textit{On the linear span of projections in certain simple C*-algebras,}
Indiana Univ. Math. J., \textbf{51} (2002) no. 3, 753--771.


\bibitem{MarcouxSmallNumberCommutators} L. W. Marcoux,
\textit{Sums of small number of commutators,} J. Operator Theory,
\textbf{56} (2006) no. 1, 111--142.


\bibitem{MarcouxIrishSurvey} L. W. Marcoux,
\textit{Projections, commutators and Lie ideals in C*-algebras,}
Math. Proc. R. Ir. Acad., \textbf{110A} (2010) no. 1, 31--55.

\bibitem{MarcouxMurphy}  L. W. Marcoux and G. J. Murphy,
\textit{Unitarily-invariant linear spaces in C*-algebras,} Proc. Amer. Math.
Soc., \textbf{126} (1998) 3597--3605.



\bibitem{OrtegaRordamThiel}  E. Ortega, M. Rordam and H. Thiel,
\textit{The Cuntz semigroup and comparison of open projections,}
(2010), Preprint,  A copy is available at
http://arxiv.org/pdf/1008.3497.

\bibitem{PearcyToppingIdempotents}  C. Pearcy and D. Topping, \textit{Sums of small numbers of idempotents,} Michigan Math. J., \textbf{14} (1967) 453--465.

\bibitem{PearcyTopping}  C. Pearcy and D. Topping, \textit{Commutators and
certain $II_1$-factors,} J. Funct. Anal., \textbf{3} (1969) 69--78.

\bibitem{PereraIdeals} F. Perera, \textit{Ideal Structure of multiplier algebras
of simple C*-algebras with real rank zero.} Canad. J. Math. \textbf{53}, (2001) no 3,  592--630.

\bibitem{Pop} C. Pop, \textit{Finite sums of commutators,}
Proc. Amer. Math. Soc., \textbf{130} (2002) 3039--3041.

\bibitem{RordamIdeals}M. Rordam, \textit{Ideals in the multiplier algebra of a stable C*-algebra.} \textbf{25}(1991)
J. Operator Theory, no. 2, 283--298.

\bibitem{ThomsenCommutator} K. Thomsen, \textit{Finite sums and products of
commutators in inductive limit C*-algebras,}  Ann. Inst. Fourier, Grenoble,
\textbf{43} (1993) no. 1, 225--249.

\bibitem{ZhangRiesz} S. Zhang, \textit{A Riesz decomposition property and
ideal structure of multiplier algebras,} J. Operator Theory, \textbf{24}
(1990) 209--226.

\bibitem{Zhang} S. Zhang, \textit{On the structure of projections and ideals
of corona algebras,} Canad. J. Math., \textbf{41} (1989) no. 4, 721-742.

\bibitem{ZhangWeylVonNeumann} S. Zhang,
 \textit{Certain C*-algebras with real rank
zero and their corona and multiplier algebras. I,} Pacific J. Math.,
\textbf{155} (1992) no. 1, 169--197.

\bibitem{ZhangMatricial} S. Zhang, \textit{Matricial structure and homotopy type of simple C*-algebras with real rank zero,} J. Operator Theory,
\textbf{26} (1991) no. 2, 283- 312.
\end{thebibliography}
\end{document}